\documentclass[a4paper,11pt,twoside,leqno]{article}

\usepackage{amsmath,amsthm,amsfonts,latexsym,amscd,amssymb}
\numberwithin{equation}{section}
\input xypic
\def\qed{{\hbadness=10000\hfill\ \vbox{\hrule height.09ex
			\hbox{\vrule width.09ex height1.55ex depth.2ex \kern1.8ex
				\vrule width.09ex height1.55ex depth.2ex}\hrule height.09ex}\break
		\bigskip}}
\newtheorem{theorem}{Theorem}[section]

\newtheorem{corollary}{Corollary}[section]

\theoremstyle{definition}

\theoremstyle{remark}

\newcommand{\n}{\noindent}

\begin{document}
	
	\linespread{1}\title {\textbf{On Relatively Normal-Slant Helices and Isophotic Curves}}
	
	\author{Akhilesh Yadav, Buddhadev Pal}
	\date{}
 
\maketitle 

	\noindent\textbf{Abstract:} In this paper, we give smoe characterizations of relatively normal-slant helices and isophotic curves on a smooth surface immersed in Euclidean 3-space with respect to their position vevtor. We also introduce the methods for generating an isophotic curve on a given surface by its  parametric or implicit equation.

\n\textbf {Mathematics Subject Classification (2010):} 53A04, 53A05.

\n\textbf{Key words:} Frenet-frame, Darboux-frame, Relatively Normal-slant helices, Isophotic curve.

\section{Introduction}
 \n In [2], Chen B. Y., introduced the notion of rectifying curve  as a space curve whose position vector always lies in its rectifying plane. He also studied several characterizations of rectifying curves, which enable us to interpret rectifying curves kinematically, as those curves whose position vector field determines the axis of instantaneous rotation at each point of the curve. In 2004, Izumiya and Takeuchi defined a slant helix in $E^3$ by the property that the principal normal vector makes a constant angle with a fixed direction and obtained a necessary and sufficient condition for a curve $\gamma$ with $\kappa(s)>0$ to be a slant helix [6]. In [1], Altunkaya and Kula studied rectifying slant helices and found the position vector of these curves.
 
 \n On the other hand, when we study space curve on a smooth surface immersed in Euclidean 3-space at every point of the curve a moving orthonormal frame called Darboux frame $\left\lbrace T, V, U\right\rbrace$ comes naturally. In [3], authors gave some characterizations of position vector of a unit speed curve in a regular surface immersed in Euclidean 3-space which always lies in the planes spanned by $\left\lbrace T, U\right\rbrace$, $\left\lbrace T, V\right\rbrace$ and $\left\lbrace V, U\right\rbrace$, respectively by using the Darboux frame. In [8], Macit N. and Duldul M., introduced the notion of relatively normal-slant helix as a curve whose vector
 field $V$ makes a constant angle with a fixed direction and gave some characterizations for such curves. In [5], Dogan F. and Yayli Y. studied isophotic curves on a surface in Euclidean $3$-space and found the axis of an isophotic curve via its Darboux frame and also gave some characterizations about the isophotic curve and its axis. 
 
 \n In [7], authors studied special vector fields along a curve associated to the Darboux frame and investigate their singularities as an application of the theory of spherical dualities. Moreover, they gave characterizations of isophotic curves on a surface by using one of the special vector fields. In this paper, we study relatively normal-slant helices and isophotic curves on a smooth surface immersed in Euclidean 3-space. The paper is arrange as follws: In section 2, we discuss some basic theory of unit speed parametrized curve on a smooth surface. In section 3, we study relatively normal-slant helices on a smooth surface with position vectors lying in the plane spanned by $\left\lbrace T, U\right\rbrace$. Section 4 is devoted to the study of isophotic curves on a smooth surface  with position vectors lying in the plane spanned by $\left\lbrace T, V\right\rbrace$. In subsection 4.1 and 4.2, we introduce some methods for generating the isophotic curve with the chosen direction and constant angle on a given surface by its  parametric and implicit equation, respectively.
 
\section{Preliminaries}
\n Let $\gamma: I \rightarrow E^3$, where $I = (\alpha, \beta)\subset\mathbb{R}$, be the unit speed parametrized curve that has at least four continuous derivatives. Then the tangent vector of the curve $\gamma$ be denoted by $T$ and given by $T(s) = \gamma^{'}(s)$, $\forall s \in I$,  where $\gamma^{'}$ denote the derivative of $ \gamma$ with respect to the arc length
parameter $s$. The binormal vector $B$ is defined by $B = T\times N$, where $N$ is the principal normal vector to the curve $\gamma$. The Frenet-Serret equations are given by
\begin{equation}
T^{'}(s) = \kappa(s)N(s),
\end{equation}
\begin{equation}
N^{'}(s) = -\kappa(s)T(s) + \tau(s)B(s),
\end{equation}
\begin{equation}
B^{'}(s) = -\tau(s)N(s),
\end{equation}
where $\kappa(s)$ and $\tau(s)$ are smooth functions of $s$, called curvature and torsion of the curve $\gamma$.

\n Let $\sigma: S\subset\mathbb{R}^2\rightarrow M$ be the coordinate chart for a smooth surface $M$ immersed in Euclidean space $E^3$ and the unit speed parametrized curve $\gamma: I \rightarrow M\subset E^3$, where $I = (\alpha, \beta)\subset\mathbb{R}$, contained in the image of a surface patch $\sigma$ in the atlas of $M$. Then $\gamma(s)$ is given by 
\begin{equation}
\gamma(s) = \sigma(u(s), v(s)),\indent \forall s \in I.
\end{equation}
Now, the curve $\gamma(s)$ lies on the surface $M$ then the Darboux frame $\left\lbrace T, V, U\right\rbrace$ at each point of the curve $\gamma(s)$ is given as follows:

\begin{equation}
\left[ {\begin{array}{c}
	T^{'} \\
	V^{'} \\
	U ^{'}\\
	\end{array} } \right]=
\left[ {\begin{array}{ccc}
	0 & k_g & k_n \\
	-k_g & 0 & \tau_{g}\\
	-k_n & -\tau_{g} & 0\\
	\end{array} } \right]
\left[ {\begin{array}{c}
	T \\
	V \\
	U \\
	\end{array} } \right],
\end{equation}
where $\kappa_{g}$, $\kappa_{n}$ and $\tau_{g}$ are the geodesic curvature, normal curvature and geodesic torsion, respectively. 

 \n Again, since $\gamma(s)$ is unit-speed curve lies on surface $M$, $\gamma^{''}$ is perpendicular to $\gamma^{'}(= T)$, and hence is a linear combination of $U$ and $V (= U \times T)$. Thus 
 \begin{equation}
 \gamma^{''}(s) = k_n(s) U(s) + k_g(s) V(s),
 \end{equation}
 As $U$ and $V$ are perpendicular unit vectors therefore from (2.6), we get
 \begin{equation}
 k_n(s) = \gamma^{''}(s).U(s)\indent and \indent k_g(s) = \gamma^{''}(s).V(s).
 \end{equation}
 Also from (2.1) and (2.7), we obtain
  \begin{equation}
 k_n(s) = \kappa(s)N(s).U(s)\indent and \indent k_g(s) = \kappa(s)N(s).V(s),
 \end{equation}
 which implies
 \begin{equation}
 k_n(s) = \kappa(s)sin\theta\indent and \indent k_g(s) = \kappa(s)cos\theta,
 \end{equation}
  where $\theta$ is the angle between vectors $N$ and $V$. Thus the curve $\gamma$ is a geodesic curve if and only if $k_g = 0$ and the curve $\gamma$ is an asymptotic line if and only if $k_n = 0$. Also geodesic torsion $\tau_g$ is given by $\tau_g = \tau - \theta^{'}$.
  
 \n Now, Differentiating (2.4) with respect to $s$, we get
 \begin{equation}
 T(s)= \gamma^{'}(s) = u^{'}\sigma_u + v^{'}\sigma_v.
 \end{equation}
 As the curve $\gamma$ is unit speed curve on $M$. Thus
  \begin{equation}
  E{u^{'}}^2 + 2Fu^{'}v^{'} + G{v^{'}}^2 = 1.
 \end{equation}
 The unit normal $U$ to the surface $M$ is given by
 \begin{equation}
 U(s) = \frac{\sigma_u\times\sigma_v}{||\sigma_u\times\sigma_v||} = \frac{\sigma_u\times\sigma_v}{\sqrt{EG-F^2}}.
 \end{equation}
 Also, since $V = U \times T$ by using (2.10) and (2.11), we obtain
  \begin{equation}
 V(s) = \frac{1}{\sqrt{EG-F^2}}(Eu^{'}\sigma_v + F(v^{'}\sigma_v - u^{'}\sigma_u) -G v^{'}\sigma_u),
 \end{equation}
 where $E = \sigma_u.\sigma_u$, $F = \sigma_u.\sigma_v$ and $G = \sigma_v.\sigma_v$ are coefficients of first fundamental form.
 
\section{Relatively normal-slant helices on a smooth surface}

\n Let $\gamma$ be a unit speed curve on an oriented surface $M$ and $(T, V, U)$ be the Darboux frame along $\gamma(s)$. The curve $\gamma$ is called a relatively normal-slant helix if the vector field $V$ of $\gamma$ makes a constant angle with a fixed direction $[8]$ , i.e. there exists a fixed unit vector $d$ and a constant angle $\phi$ such that $<V, d> = cos\phi$. The unit vector $d$ is called the axis of the relatively normal-slant helix. In this section, we study relatively normal-slant helices on a smooth surface immersed in Euclidean $3$-space, whose position vector lies in the plane spanned by $\left\lbrace T, U\right\rbrace$.

\n If position vector  of the curve $\gamma$ always lies in the plane $sp\left\lbrace T, U\right\rbrace$ then the position vector of the curve satisfies the equation 
\begin{equation}
\gamma(s) = \lambda_1(s)T(s) + \lambda_2(s)U(s),
\end{equation}
for some differentiable functions $\lambda_1(s)$ and $\lambda_2(s)$. Then from $[3]$, we have
\begin{equation}
\lambda_1(s) = \frac{c\tau_{g}}{\kappa_{g}}e^{-\int\frac{\tau_{g}\kappa_{n}}{\kappa_{g}}ds},
\end{equation}
\begin{equation}
\lambda_2(s) = ce^{-\int\frac{\tau_{g}\kappa_{n}}{\kappa_{g}}ds},
\end{equation}
where $c \in\mathbb{R}_0$.
Also from $[3]$, the curvature functions $\kappa_{n}(s)$, $\kappa_{g}(s)$ and $\tau_{g}(s)$ satisfy the equation
\begin{equation}
\left( \left( \dfrac{\tau_g}{k_g}\right) ^{'} - \left( \left( \dfrac{\tau_g}{k_g}\right) ^2 + 1 \right) k_n\right)(s) = \frac{1}{c}e^{\int\frac{\tau_{g}\kappa_{n}}{\kappa_{g}}ds}.
\end{equation}
\begin {theorem}$[8]$ A unit speed curve $\gamma$ on a surface $M$ with $(\tau_g(s), k_g(s)) \neq (0, 0)$ is a relatively normal-slant helix if and only if
\begin{equation}
\mu_v(s) = (\dfrac{1}{(k_g^2 + \tau_g^2)^{\frac{3}{2}}}(k_g\tau_g^{'} - \tau_gk_g^{'} - k_n(k_g^2 + \tau_g^2)))(s)
\end{equation}
is a constant function.
\end {theorem}
\begin {theorem}  Let $\gamma$ be a unit speed curve on a smooth surface $M$ with position vector lies in the plane $sp\left\lbrace T, U\right\rbrace$. Then the curve $\gamma$ is a relatively normal-slant helix if and only if $\dfrac{k_g^2}{(k_g^2 + \tau_g^2)^{\frac{3}{2}}}e^{\int\frac{\tau_{g}\kappa_{n}}{\kappa_{g}}ds}$ is a constant function.
\end{theorem}
\n\begin{proof} First assume that curve $\gamma$ is a relatively normal-slant helix then 
\begin{equation}
\dfrac{k_g^2}{(k_g^2 + \tau_g^2)^{\frac{3}{2}}}\left( \left( \dfrac{\tau_g}{k_g}\right) ^{'} - \left( \left( \dfrac{\tau_g}{k_g}\right) ^2 + 1 \right) k_n\right)(s) 
\end{equation}
is a constant function.
Then from (3.4) and (3.6), we get $\dfrac{k_g^2}{(k_g^2 + \tau_g^2)^{\frac{3}{2}}}e^{\int\frac{\tau_{g}\kappa_{n}}{\kappa_{g}}ds}$ is a constant function.

\n Conversely, suppose that $\dfrac{k_g^2}{(k_g^2 + \tau_g^2)^{\frac{3}{2}}}e^{\int\frac{\tau_{g}\kappa_{n}}{\kappa_{g}}ds}$ is a constant function. Then by using this in (3.4), we get the required result.
\end{proof}
\begin {corollary}  Let $\gamma$ be a unit speed asymptotic curve on a smooth surface $M$ with position vector lies in the plane $sp\left\lbrace T, U\right\rbrace$. Then the curve $\gamma$ is a relatively normal-slant helix if and only if $\left( \dfrac{k^2}{(k^2 + \tau^2)^{\frac{3}{2}}}\right)(s)$ is a constant function.
\end{corollary}
\begin {corollary}  Let $\gamma$ be a unit speed line of curvature on a smooth surface $M$ with position vector lies in the plane $sp\left\lbrace T, U\right\rbrace$. Then the curve $\gamma$ is a relatively normal-slant helix if and only if $\kappa_{g}(s)$ is a constant function.
\end{corollary}
\begin {theorem}  Let $\gamma$ be a relatively normal-slant helix  on a smooth surface $M$ with position vector lies in the plane $sp\left\lbrace T, U\right\rbrace$. Then the position vector of the curve $\gamma$ satisfies the equation 
\begin{equation}
\gamma(s) = \dfrac{k_g\tau_g}{(k_g^2 + \tau_g^2)^{\frac{3}{2}}}T(s) + \dfrac{k_g^2}{(k_g^2 + \tau_g^2)^{\frac{3}{2}}}U(s).
\end{equation}
\end{theorem}
\n \begin{proof} Since position vector of $\gamma$ lies in the plane $sp\left\lbrace T, U\right\rbrace$ therefore from (3.1), (3.2) and (3.3), we have
\begin{equation}
\gamma(s) = \frac{c\tau_{g}}{\kappa_{g}}e^{-\int\frac{\tau_{g}\kappa_{n}}{\kappa_{g}}ds}T(s) + ce^{-\int\frac{\tau_{g}\kappa_{n}}{\kappa_{g}}ds}U(s).
\end{equation}
Also, $\gamma$ is relatively normal-slant helix with position vector in the plane $sp\left\lbrace T, U\right\rbrace$. Then, we have 
\begin{equation} 
\dfrac{k_g^2}{(k_g^2 + \tau_g^2)^{\frac{3}{2}}}e^{\int\frac{\tau_{g}\kappa_{n}}{\kappa_{g}}ds} = c (const.).
\end{equation} 
\n Thus from (3.8) and (3.9), we obtain the required result.
\end{proof}

\begin {corollary}  Let $\gamma$ be a relatively normal-slant helix  on a smooth surface $M$ with position vector lies in the plane $sp\left\lbrace T, U\right\rbrace$. If $\gamma$ is asymptotic curve on $M$ then $<\gamma, T> = \dfrac{k\tau}{(k^2 + \tau^2)^{\frac{3}{2}}}$ and $<\gamma, U> =  \dfrac{k^2}{(k^2 + \tau^2)^{\frac{3}{2}}}$.
\end{corollary}
\begin {corollary}  Let $\gamma$ be a relatively normal-slant helix  on a smooth surface $M$ with position vector lies in the plane $sp\left\lbrace T, U\right\rbrace$. If $\gamma$ is line of curvature on $M$ then $<\gamma, T> = 0$ and $<\gamma, U> = \frac{1}{k_g}$.
\end{corollary}
\begin {theorem}  Let $\gamma$ be a relatively normal-slant helix with $k_n = 0$ and position vector of $\gamma$ lies in the plane $sp\left\lbrace T, U\right\rbrace$. Then $\gamma$ is slant helix if and only if $\gamma$ is rectifying curve.
\end{theorem}
\n \begin{proof} Let $\gamma$ be a relatively normal-slant helix with $k_n = 0$ and position vector of $\gamma$ lies in the plane $sp\left\lbrace T, U\right\rbrace$. Suppose $\gamma$ is slant helix then from [6], we have
\begin{equation}
\mu(s) =\frac{\kappa^{2}}{(\kappa^2 + \tau^2)^{3/2}}(\frac{\tau}{\kappa})^{'}(s)= constant.
\end{equation}
Now, from corollary 3.1 and equation (3.10), we get $(\frac{\tau}{\kappa})^{'}(s)= constant$, which implies $(\frac{\tau}{\kappa})(s)= c_1s + c_2$, where $c_1, c_2 \in\mathbb{R}_0$. Thus $\gamma$ is a rectifying curve.
Conversely suppose that $\gamma$ is a rectifying curve. Then $(\frac{\tau}{\kappa})^{'}(s)= constant$. Also since $\gamma$ is relatively normal-slant helix with $k_n = 0$ and position vector of $\gamma$ lies in the plane $sp\left\lbrace T, U\right\rbrace$ therefore $\left( \dfrac{k^2}{(k^2 + \tau^2)^{\frac{3}{2}}}\right)(s)$ is a constant function. Thus $\left( \dfrac{k^2}{(k^2 + \tau^2)^{\frac{3}{2}}}\right)(\frac{\tau}{\kappa})^{'}(s)$ is a constant function. Hence $\gamma$ is slant helix.
\end{proof}

\section{Isophotic curves on a smooth surface}

\n  Let $\gamma$ be a unit speed curve on an oriented surface $M$. The curve $\gamma$ is called an isophotic curve if the unit normal vector field $U$ of $M$ along $\gamma$ makes a constant angle with a fixed direction $[5]$, i.e. there exists a fixed unit vector $d$ and a constant angle $\phi$ such that $<U, d> = cos\phi$. The unit vector $d$ is called the axis of the isophotic curve. In this section, we study isophotic curves on a smooth surface immersed in Euclidean $3$-space, whose position vector lies in the plane spanned by $\left\lbrace T, V\right\rbrace$. 

\n If position vector  of the curve $\gamma$ always lies in the plane $sp\left\lbrace T, U\right\rbrace$ then the position vector of the curve satisfies the equation 
\begin{equation}
\gamma(s) = \mu_1(s)T(s) + \mu_2(s)V(s),
\end{equation}
for some differentiable functions $\mu_1(s)$ and $\mu_2(s)$. Then from $[3]$, we have
\begin{equation}
\mu_1(s) = \frac{-c\tau_{g}}{\kappa_{n}}e^{\int\frac{\tau_{g}\kappa_{g}}{\kappa_{n}}ds},
\end{equation}
\begin{equation}
\mu_2(s) = ce^{\int\frac{\tau_{g}\kappa_{g}}{\kappa_{n}}ds},
\end{equation}
where $c \in\mathbb{R}_0$.
Also from $[3]$, the curvature functions $\kappa_{n}(s)$, $\kappa_{g}(s)$ and $\tau_{g}(s)$ satisfy the equation
\begin{equation}
\left( \left( \dfrac{\tau_g}{k_n}\right) ^{'} - \left( \left( \dfrac{\tau_g}{k_n}\right) ^2 + 1 \right) k_g\right)(s) = -\frac{1}{c}e^{-\int\frac{\tau_{g}\kappa_{g}}{\kappa_{n}}ds}.
\end{equation}
\begin {theorem}$[5]$ A unit speed curve $\gamma$ on a surface $M$ with $(\tau_g(s), k_n(s)) \neq (0, 0)$ is an isophotic curve if and only if
\begin{equation}
\mu_u(s) = (\dfrac{1}{(k_n^2 + \tau_g^2)^{\frac{3}{2}}}(k_n\tau_g^{'} - \tau_gk_n^{'} - k_g(k_n^2 + \tau_g^2)))(s)
\end{equation}
is a constant function.
\end {theorem}
\begin {theorem}  Let $\gamma$ be a unit speed curve on a smooth surface $M$ with position vector lies in the plane $sp\left\lbrace T, V\right\rbrace$. Then the curve $\gamma$ is an isophotic curve if and only if $\dfrac{k_n^2}{(k_n^2 + \tau_g^2)^{\frac{3}{2}}}e^{-\int\frac{\tau_{g}\kappa_{g}}{\kappa_{n}}ds}$ is a constant function.
\end{theorem}
\n \begin{proof} First assume that curve $\gamma$ is an isophotic curve, then 
\begin{equation}
\dfrac{k_n^2}{(k_n^2 + \tau_g^2)^{\frac{3}{2}}}\left( \left( \dfrac{\tau_g}{k_n}\right) ^{'} - \left( \left( \dfrac{\tau_g}{k_n}\right) ^2 + 1 \right) k_g\right)(s) 
\end{equation}
is a constant function.
Then from (4.4) and (4.6), we get $\dfrac{k_n^2}{(k_n^2 + \tau_g^2)^{\frac{3}{2}}}e^{-\int\frac{\tau_{g}\kappa_{g}}{\kappa_{n}}ds}$ is a constant function.

\n Conversely, suppose that $\dfrac{k_n^2}{(k_n^2 + \tau_g^2)^{\frac{3}{2}}}e^{-\int\frac{\tau_{g}\kappa_{g}}{\kappa_{n}}ds}$ is a constant function. Then by using this in (4.4), we obtain (4.6). Hence $\gamma$ is an isophotic curve.
 \end{proof}
\begin {corollary}  Let $\gamma$ be a unit speed geodesic curve on a smooth surface $M$ with position vector lies in the plane $sp\left\lbrace T, V\right\rbrace$. Then the curve $\gamma$ is an isophotic curve if and only if $\left( \dfrac{k^2}{(k^2 + \tau^2)^{\frac{3}{2}}}\right)(s)$ is a constant function.
\end{corollary}
\begin {corollary}  Let $\gamma$ be a unit speed line of curvature on a smooth surface $M$ with position vector lies in the plane $sp\left\lbrace T, V\right\rbrace$. Then the curve $\gamma$ is an isophotic curve if and only if $\kappa_{n}(s)$ is a constant function.
\end{corollary}
\begin {theorem}  Let $\gamma$ be an isophotic curve  on a smooth surface $M$ with position vector lies in the plane $sp\left\lbrace T, V\right\rbrace$. Then the position vector of the curve $\gamma$ satisfies the equation 
\begin{equation}
\gamma(s) = \dfrac{k_n^2}{(k_n^2 + \tau_g^2)^{\frac{3}{2}}}V(s) - \dfrac{k_n\tau_g}{(k_n^2 + \tau_g^2)^{\frac{3}{2}}}T(s).
\end{equation}
\end{theorem}
\n \begin{proof} Since position vector of $\gamma$ lies in the plane $\left\lbrace T, V\right\rbrace$ therefore from (4.1), (4.2) and (4.3), we have
\begin{equation}
\gamma(s) = -\frac{c\tau_{g}}{\kappa_{n}}e^{\int\frac{\tau_{g}\kappa_{g}}{\kappa_{n}}ds}T(s) + ce^{\int\frac{\tau_{g}\kappa_{g}}{\kappa_{n}}ds}V(s).
\end{equation}
Also, $\gamma$ is an isophotic curve with position vector in the plane $sp\left\lbrace T, V\right\rbrace$. Then, we have 
\begin{equation} 
\dfrac{k_n^2}{(k_n^2 + \tau_g^2)^{\frac{3}{2}}}e^{-\int\frac{\tau_{g}\kappa_{g}}{\kappa_{n}}ds} = c (const.).
\end{equation} 
\n Thus from (4.8) and (4.9), we obtain the required result.
 \end{proof}

\begin {corollary}  Let $\gamma$ be an isophotic curve on a smooth surface $M$ with position vector lies in the plane $sp\left\lbrace T, V\right\rbrace$. If $\gamma$ is geodesic curve on $M$ then $<\gamma, T> = -\dfrac{k\tau}{(k^2 + \tau^2)^{\frac{3}{2}}}$ and $<\gamma, V> =  \dfrac{k^2}{(k^2 + \tau^2)^{\frac{3}{2}}}$.
\end{corollary}
\begin {corollary}  Let $\gamma$ be an isophotic curve on a smooth surface $M$ with position vector lies in the plane $sp\left\lbrace T, V\right\rbrace$. If $\gamma$ is line of curvature on $M$ then $<\gamma, T> = 0$ and $<\gamma, V> = \frac{1}{k_n}$.
\end{corollary}
\begin {theorem}  Let $\gamma$ be an isophotic curve with $k_g = 0$ and position vector of $\gamma$ lies in the plane $sp\left\lbrace T, V\right\rbrace$. Then $\gamma$ is slant helix if and only if $\gamma$ is rectifying curve.
\end{theorem}
\n \begin{proof} Let $\gamma$ be an isophotic curve with $k_g = 0$ and position vector of $\gamma$ lies in the plane $sp\left\lbrace T, V\right\rbrace$. Suppose $\gamma$ is slant helix then from [6], we have $\mu(s) =\frac{\kappa^{2}}{(\kappa^2 + \tau^2)^{3/2}}(\frac{\tau}{\kappa})^{'}(s)= constant$.
Thus, by using corollary 4.1 in above equation, we get $(\frac{\tau}{\kappa})^{'}(s)= constant$, which implies $(\frac{\tau}{\kappa})(s)= c_1s + c_2$, where $c_1, c_2 \in\mathbb{R}_0$. Hence, $\gamma$ is a rectifying curve.

\n Conversely suppose that $\gamma$ is a rectifying curve. Then $(\frac{\tau}{\kappa})^{'}(s)= constant$. Also since $\gamma$ is an isophotic curve with $k_g = 0$ and position vector of $\gamma$ lies in the plane $sp\left\lbrace T, V\right\rbrace$ therefore $\left( \dfrac{k^2}{(k^2 + \tau^2)^{\frac{3}{2}}}\right)(s)$ is a constant function. Thus $\left( \dfrac{k^2}{(k^2 + \tau^2)^{\frac{3}{2}}}\right)(\frac{\tau}{\kappa})^{'}(s)$ is a constant function. Hence $\gamma$ is slant helix.
 \end{proof}
\subsection{\textbf{Isophotic curves on a parametric surface}}
Let $M$ be a regular oriented surface in $E^3$ with the parametrization $\sigma = \sigma(u, v)$. Now our goal is to give the method which enables us to find the isophotic curve $\gamma(s) = \sigma(u(s), v(s))$ (if exists) lying on $M$, when axis $d$ of the isophotic curve and the constant angle $\phi$ are given.

\n Then we have $<U, d> = cos\phi$, which implies $<U^{'}, d> = 0$. By using (2.5) in the above equation, we get
\begin{equation} 
<k_nT +\tau_gV, d> = 0.
\end{equation} 
Then by using (2.10) and (2.13) in (4.10), we obtain
\begin{equation} \begin{split}
&\left(\sqrt{EG-F^2}k_n<\sigma_u, d> + E\tau_g<\sigma_v, d> - F\tau_g<\sigma_u, d>\right)\dfrac{du}{ds} \\&+ \left(\sqrt{EG-F^2}k_n<\sigma_v, d> + F\tau_g<\sigma_v, d> -G\tau_g<\sigma_u, d>\right)\dfrac{dv}{ds}  = 0.
\end{split}\end{equation} 
Now, from (2.11) and (4.11), we get
\begin{equation} 
\dfrac{du}{ds}  =  \pm\dfrac{{\bigtriangleup^*}^2}{({\bigtriangleup^*}^2E +-2{\bigtriangleup^*}{\bigtriangleup}F + {\bigtriangleup}^2G)^\frac{1}{2}}
\end{equation} 
and
\begin{equation} 
\dfrac{dv}{ds}  =  \pm\dfrac{{\bigtriangleup}^2}{({\bigtriangleup^*}^2E +-2{\bigtriangleup^*}{\bigtriangleup}F + {\bigtriangleup}^2G)^\frac{1}{2}},
\end{equation} 
where $$\bigtriangleup = \sqrt{EG-F^2}k_n<\sigma_u, d> + E\tau_g<\sigma_v, d> - F\tau_g<\sigma_u, d>,$$
$$\bigtriangleup^* = \sqrt{EG-F^2}k_n<\sigma_v, d> + F\tau_g<\sigma_v, d> -G\tau_g<\sigma_u, d>.$$
If we solve the system of ODE (4.12) and (4.13) together with the initial point $u(0) = u_0$, $v(0) = v_0$, we obtain the desired isophotic curves on $M$ by substituting $u(s)$, $v(s)$ into $\gamma(s) = \sigma(u(s), v(s))$.

\n Remark: (i) If $({\bigtriangleup^*}^2E +-2{\bigtriangleup^*}{\bigtriangleup}F + {\bigtriangleup}^2G) \leq 0$, then there does not exist a isophotic curve on $M$ with the given axis and angle.

(ii) If $({\bigtriangleup^*}^2E +-2{\bigtriangleup^*}{\bigtriangleup}F + {\bigtriangleup}^2G) > 0$, then we have two isophotic curve on the surface $M$.

\subsection{\textbf{Isophotic curves lying on an implicit surface}}
Let $M$ be a surface given in implicit form by $f(x, y, z) = 0$. Let $\gamma(s) = (x(s), y(s), z(s))$ be an isophotic curve on $M$, which makes the given constant angle $\phi$ with the given axis $d = (d_1, d_2, d_3)$ and $(T, V, U)$ be its Darboux frame field. We need to find $x(s)$, $y(s)$, $z(s)$ to obtain $\gamma(s)$.

\n Now, $\gamma(s)$ is the unit speed isophotic curve on $M$. Thus
\begin{equation} 
f_x\dfrac{dx}{ds} + f_y\dfrac{dy}{ds} + f_z\dfrac{dz}{ds} = 0
\end{equation} 
and
\begin{equation} 
{\dfrac{dx}{ds}}^2 + {\dfrac{dy}{ds}}^2 + {\dfrac{dz}{ds}}^2 =1,
\end{equation} 
where $f_\xi = \dfrac{\partial f}{\partial \xi}$.

\n Also, the unit normal vector field of the surface is $U = \dfrac{\bigtriangledown f}{||\bigtriangledown f||}$. Thus
$$V = \dfrac{\bigtriangledown f}{||\bigtriangledown f||} \times T = \dfrac{1}{||\bigtriangledown f||}\left( f_y\dfrac{dz}{ds} - f_z\dfrac{dy}{ds},  f_z\dfrac{dx}{ds} - f_x\dfrac{dz}{ds}, f_x\dfrac{dy}{ds} - f_y\dfrac{dx}{ds}\right).$$
Now by putting the value of $T$ and $V$ in (4.10), we get
\begin{equation} \begin{split}
\left( d_1k_n + (d_2f_z - d_3f_y)\tau_g\right) \dfrac{dx}{ds}& + \left( d_2k_n + (d_3f_x - d_1f_z)\tau_g\right) \dfrac{dy}{ds} \\&+ \left( d_3k_n + (d_1f_y - d_2f_x)\tau_g\right) \dfrac{dz}{ds} = 0.
\end{split}\end{equation} 
From (4.14) and (4.16), we obtain
\begin{equation} 
\dfrac{dx}{ds} = \left( \dfrac{f_y\Omega_3 - f_z\Omega_2}{f_x\Omega_2 - f_y\Omega_1}\right) \dfrac{dz}{ds}
\end{equation}
and
\begin{equation} 
\dfrac{dy}{ds} = \left( \dfrac{f_z\Omega_1 - f_x\Omega_3}{f_x\Omega_2 - f_y\Omega_1}\right) \dfrac{dz}{ds},
\end{equation}
where $\Omega_1 = d_1k_n + (d_2f_z - d_3f_y)\tau_g$, $\Omega_2 = d_2k_n + (d_3f_x - d_1f_z)\tau_g$ and $\Omega_3 = d_3k_n + (d_1f_y - d_2f_x)\tau_g$.
Also from (4.15), (4.17) and (4.18), we get
\begin{equation} 
\dfrac{dz}{ds} = \pm\dfrac{f_x\Omega_2 - f_y\Omega_1}{\left( (f_x\Omega_2 - f_y\Omega_1)^2 + (f_y\Omega_3 - f_z\Omega_2)^2 + (f_z\Omega_1 - f_x\Omega_3)^2\right) \frac{1}{2}}.
\end{equation}
If we substitute (4.19) into (4.17) and (4.18), we obtain an explicit $1^{st}$ order ordinary differential equation system. Thus, together with the initial point $x(0) = x_0$, $y(0) = y_0$, $z(0) = z_0$, we have an initial value problem. The solution of this problem gives isophotic curve on $M$.

Remark. (i) If $(f_x\Omega_2 - f_y\Omega_1)^2 + (f_y\Omega_3 - f_z\Omega_2)^2 + (f_z\Omega_1 - f_x\Omega_3)^2 = 0$ at the point $(x_0, y_0, z_0)$, then there does not exist any isophotic curve with the given direction $d$ and angle $\phi$.

(ii) If $(f_x\Omega_2 - f_y\Omega_1)^2 + (f_y\Omega_3 - f_z\Omega_2)^2 + (f_z\Omega_1 - f_x\Omega_3)^2 \neq 0$ at the point $(x_0, y_0, z_0)$, then we have two isophotic curves passing through the initial point.

\noindent\author{Akhilesh Yadav}\\
\date{Department of Mathematics, Institute of Science, \\Banaras Hindu University, Varanasi-221005, India}\\
\maketitle {\noindent E-mails: akhilesha68@gmail.com,     akhilesh$\_$mathau@rediffmail.com} 

\noindent\author{Buddhadev Pal}\\
\date{Department of Mathematics, Institute of Science, \\Banaras Hindu University, Varanasi-221005, India}\\
\maketitle {\noindent E-mail: pal.buddha@gmail.com}

\end{document}